\documentclass{amsart}
\usepackage{amssymb}

\newcommand{\KZ}{{\mathit{KZ}}}
\newcommand{\CM}{{\mathit{CM}}}
\newcommand{\R}{{\mathbb{R}}}
\newcommand{\C}{{\mathbb{C}}}

\newcommand{\Z}{{\mathbb{Z}}}
\newcommand{\dontprint}[1]{\relax}
\newtheorem%
{thm}{Theorem}[section]
\newtheorem%
{proposition}[thm]{Proposition}
\newtheorem%
{lemma}[thm]{Lemma}
\newtheorem%
{lemmadef}[thm]{Lemma-Definition}
\newtheorem%
{corollary}[thm]{Corollary}
\newtheorem%
{conjecture}[thm]{Conjecture}

\begin{document}

\title
[Action of Coxeter groups on $m$-harmonic polynomials]
{Action of Coxeter groups on $m$-harmonic polynomials and KZ equations}
%
%
       \author{G. Felder}
\address{Department of Mathematics,
          ETH, Zurich, Switzerland}
       \email{felder@math.ethz.ch}
         \author{A.P. Veselov}
\address{Department of Mathematical Sciences,
          Loughborough University, Loughborough,
          Leicestershire, LE11 3TU, UK
          and Landau Institute for Theoretical Physics, Moscow, Russia}
         \email{A.P.Veselov@lboro.ac.uk}

\begin{abstract}
The Matsuo--Cherednik correspondence is an isomorphism
from solutions of Kni\-zhnik--Zamolo\-dchi\-kov equations to eigenfunctions of 
generalized Calogero--Moser systems associated to
Coxeter groups $G$ and a multiplicity function $m$ on their
root systems. 
We apply a version of this correspondence to the most degenerate case of
zero spectral parameters. The space of eigenfunctions
is then the space $H_m$ 
of $m$-harmonic polynomials, recently introduced
in \cite{{FeiginVeselov}}.
 We compute the Poincar\'e polynomials for the
space $H_m$ and for its isotypical components
corresponding to each irreducible
representation of the group $G$. We also
give an explicit formula
for $m$-harmonic polynomials of lowest positive degree in the
$S_n$ case.
\end{abstract}

\dontprint{ 
The Matsuo-Cherednik correspondence is an isomorphism from solutions
of Knizhnik-Zamolodchikov equations to eigenfunctions of generalized
Calogero-Moser systems associated to Coxeter groups G and a
multiplicity function m on their root systems. 
We apply a version of this correspondence to the most degenerate case of
zero spectral parameters. The
space of eigenfunctions is then the space H_m of m-harmonic
polynomials, recently introduced in math-ph/0105014. We compute the
Poincare' polynomials for the space H_m and of its isotypical
components corresponding to each irreducible representation of the
group G. We also give an explicit formula for m-harmonic polynomials
of lowest positive degree in the S_n case.
}

\maketitle

\section{Introduction}
Let $G$ be a finite group generated by orthogonal hyperplane
reflections of  $n$-dimensional Euclidean space $\R^n$.
We label all reflections
$s_\alpha$ by a set of arbitrarily chosen normal vectors $\alpha\in A$.
We also fix a system of multiplicities $m_\alpha\in\Z_{\geq0}$ such that
$m_\alpha=m_\beta$ if $s_\alpha$ and $s_\beta$ belong to the same
conjugacy class.
The group $G$ also acts on the complexification $V=\C^n$ of $\R^n$.
We denote
by $\Pi_\alpha$ the complexification of the reflection hyperplane of
$s_\alpha$.

By Chevalley's theorem, the algebra
of invariants $S^G$ in the graded algebra $S=S(V)=\oplus_{j=0}^\infty S^j$ of
polynomial functions on $V$ is freely generated by homogeneous polynomials
$\sigma_1,\dots,\sigma_n$, $\sigma_i\in S^{d_i}$.
We identify $V$ with its dual space by means of
the bilinear form $(\ ,\ )$ induced by the Euclidean inner product.
So we identify the symmetric algebra $S$
both with polynomial
functions $z\mapsto p(z)$ on $V$ and with differential operators
$p(\partial)$ with constant coefficients. Under this identification,
a vector $\xi\in V$ corresponds to the linear function $z\mapsto(\xi,z)$
and to the derivative $\partial_\xi$ in the direction of $\xi$.

Let us consider the following system of equations
\begin{equation}\label{0-har}
\sigma_i(\partial)\phi= 0, \qquad i=1,\dots,n,
\end{equation}
where $\sigma_1,\dots,\sigma_n$ as above are the generators of $S^G$.
The solutions of this system form a $|G|$-dimensional space $H$ and
turn out to be polynomial
\cite{Steinberg2}.
These polynomials are called {\it harmonic}
and play an important role in the theory of Coxeter groups and symmetric spaces
(see \cite{Helgason}). It is known that they freely generate $S(V)$ as a module
over $S^G$ and the corresponding Poincar\'e polynomial $P(H,t) = \sum_k
t^k\,{\mathrm{dim}\, H^k} ,$
where $H^k$ is the subspace of harmonic polynomials of degree $k$, has the form
$$P(H,t) = \prod_{k=1}^n\frac{1-t^{d_k}}{1-t},$$
$d_1,...,d_n$ are the degrees of the basic invariants
$\sigma_1,\dots,\sigma_n.$

In the present paper, we investigate a generalization of harmonic
polynomials --- the so-called {\it $m$-harmonic} polynomials
introduced recently in \cite{FeiginVeselov}.

To define them let us assign to each reflection hyperplane
$\Pi_\alpha$ a nonnegative
integer $m_\alpha\in\Z_{\geq0}$
such that if $\Pi_\alpha =g(\Pi_\beta)$, $g\in G$ then $m_\alpha =m_\beta$.
Such a function $m$ on the set of the reflections is called {\it multiplicity}.
Since this function is constant on each orbit of the action of $G$ on its
reflection hyperplanes,
the number $q$ of parameters is equal to the number of such
orbits (or equivalently,
the number of conjugacy classes $C_1,..., C_q$ of reflections in $G$).
For irreducible groups $G$, $q$ is actually either 1 or 2.

Corresponding $m$-harmonic polynomials are defined as the solutions
of the following system generalizing (\ref{0-har}):
\begin{equation}\label{m-har}
\mathcal{L}_i\phi= 0,\qquad i=1,\dots,n,
\end{equation}
where $$\mathcal{L}_1 = \Delta -\sum_{\alpha\in
A}\frac{2m_\alpha}{(\alpha,x)}\partial_\alpha$$
is
(up to a gauge transformation)
the  generalized Calogero--Moser operator and $\mathcal{L}_i$ are
its quantum integrals with the highest terms $\sigma_i(\partial)$
(see the next section for details).
As it was shown in \cite{FeiginVeselov} all the solutions of this
system are polynomial
and form a $|G|$-dimensional space which is denoted as $H_m$. When
the multiplicity is zero
we have the space $H_0$ of usual harmonic polynomials.
For dihedral groups and constant multiplicity functions all
$m$-harmonic polynomials
have been described in \cite{FeiginVeselov} but for the general groups
even the question of the degrees of such polynomials was open.

Our main idea to attack this problem  goes back 
to Matsuo's and Cherednik's observations that
the Calogero--Moser system
\begin{equation}\label{CMk}
\mathcal{L}_i\phi= \sigma_i(\lambda)\phi,\qquad i=1,\dots,n.
\end{equation}
for generic $\lambda \in V$, namely if
$\prod_{\alpha\in A}(\lambda ,\alpha)
\neq 0$, is equivalent to a certain version
of the Knizhnik--Zamolodchikov (KZ) equations (see
\cite{Matsuo},\cite{Cherednik2} and Section \ref{s-KZ}).
These equations form a compatible system of first order
equations for a function $u$ on $V$ taking values in the
group algebra $\C[G]$. They have the form
\begin{equation}\label{e-KZ00}
\partial_\xi u(x)=\sum_{\alpha\in
A}m_\alpha\frac{(\alpha,\xi)}{(\alpha,x)}
(s_\alpha+1)u(x)+\pi(\xi,\lambda)u(x),\qquad \xi\in V,
\end{equation}
where $\pi(\xi,\lambda)$ is the linear endomorphism
of $\C[G]$ so that $\pi(\xi,\lambda)g=(\xi,g\lambda)g$, for
$g\in G$. Then, for generic $\lambda$, the composition
with the alternating representation  $\mu:u\mapsto 
\epsilon\circ u$ 
is an isomorphism from the sheaf of local solutions of 
the KZ equations to the sheaf of local solutions of 
the Calogero--Moser system \eqref{CMk}.

For non-generic $\lambda$, in particular for $\lambda=0$ this
construction does not work, as the map $\mu$ fails
to be an isomorphism if $(\lambda,\alpha)=0$ for some
$\alpha\in A$.

For this reason, we use a modification,
due essentially to Cherednik \cite{Cherednik2}, \cite{Cheredniklectures},
of the
KZ equations and of the isomorphism $\mu$,
which works for all $\lambda \in V$ including $\lambda = 0$.
The idea is that, instead of $\C[G]$, it is more natural
to consider the KZ equations for a function with values
in the $G$-module $S(V)/I(\lambda)$, where $I(\lambda)$ is the
ideal generated by $G$-invariant polynomials vanishing at $\lambda$. This module is (non-uniquely) isomorphic to $\C[G]$ 
for  all $\lambda\in V$. 
The KZ equations with values in $S(V)/I(\lambda)$
 still have the form \eqref{e-KZ00} with a 
suitable $\pi(\xi,\lambda)$, namely the multiplication by
$\xi\in V\subset S(V)$. These equations come with a 
 {\em 
Matsuo-Cherednik map} $\mu$ from local solutions 
to local solutions of the Calogero--Moser system
\eqref{CMk}, which is an isomorphism
for {\em all} $\lambda\in V$, see Theorem \ref{t-1}. 
The equations \eqref{e-KZ} 
are
recovered from 
the KZ equations with values in $S(V)/I(\lambda)$
via the map $S(V)/I(\lambda)\to \C[G]$
induced from the morphism sending a polynomial $p\in S(V)$ to 
$\sum_{g\in G} p(g\lambda)g$. This map is an isomorphism only for
generic $\lambda$.

This 
explicit form of the Calogero--Moser system (\ref{m-har}) as
a holonomic system seems to be
quite convenient. We show that it allows
one
not only to find the
degrees of $m$-harmonic polynomials
but to compute them explicitly in some special cases.

An application of our construction
 is a formula for the Poincar\'e polynomial $P(H_m,t)$
for the general Coxeter group $G$ and multiplicity function $m$.

Let us describe this formula, which is given as a sum of the
Poincar\'e polynomials for each isotypical
component.
Let $V_1,...,V_p$
be a list of all
inequivalent
irreducible representations of $G$. For any
representation $V_j$ we can
ask what is the multiplicity $p_k(V_j)$ of this representation in the
$G$-module of (usual) harmonic polynomials
of degree $k$ and define the corresponding Poincar\'e polynomial
$$P_j(t) = \sum_k p_k(V_j) t^k.$$
These polynomials are known for all Coxeter groups (see
\ref{ss-33} below).
They obey Poincar\'e duality
\begin{equation}\label{poinc}
P_{j^*}(t)=t^NP_j(t^{-1}),
\end{equation}
where $V_{j^*}$ is the tensor product of  $V_j$ by
the alternating representation
and $N$ is the total number of the reflections in $G$.

Let us define now for any conjugacy class of reflections $C_a$ the number
\begin{equation}\label{e-defd}
d^{-}_a(V_j)=\frac{2N_a\mathrm{dim}(V_{j,\alpha}^{-})}
{\mathrm{dim}(V_j)},\qquad
a=1,\dots,q,
\end{equation}
where $N_a$ is the number of elements in $C_a$ and $(V_{j,\alpha}^{-})$ is the
$-1$-eigenspace of the action of any $s_\alpha \in C_a$
on $V_j$.
This number can be expressed in terms of the polynomials
$P_j(t)$ in the following way. Let $P_{j\otimes a}(t)$ be the
Poincar\'e polynomial
corresponding to the representation $V_j \otimes \chi_a$
where $\chi_a$ is the one-dimensional representation of $G$
such that
$$
\chi_a(s)=\left\{\begin{array}{rl}
-1,& \text{if $s\in C_a$,}\\
1,& \text{if $s\in C_b$, $b\neq a$,}
\end{array}
\right.
$$
We show that
$$d^-_a(V_j)=N_a+
\left.\frac{d}{dt}\right|_{t=1}
\ln\frac{P_j(t)}{P_{j\otimes a}(t)},\qquad
a=1,\dots, q.$$
If the group $G$ acts transitively on its reflection hyperplanes (i.e. $q=1$)
the formula simplifies to Solomon's formula \cite{Solomon}
$$
d^-(V_j)=d^-_1(V_j)=2\left.\frac{d}{dt}\right|_{t=1}\ln{P_j(t)}.
$$

\begin{thm}\label{t-0}
The
Poincar\'e polynomial of the graded space of $m$-harmonic polynomials 
has a form
\[
P(H_m,t)=\sum_{j=1}^{p}
\mathrm{dim}(V_j)\,
t^{d^{-}_j(m)}
P_{j}(t),
\]
where $d^{-}_j(m) = \sum_{a=1}^qm_ad^-_a(V_j).$
It has degree $M = \sum_{\alpha \in A} (2 m_\alpha +1)= 
\sum_{a=1}^q N_a (2m_a +1)$ and
satisfies the palindromic relation
$P(H_m,t) = t^M P(H_m, t^{-1})$.
\end{thm}

For example, in the case of the symmetric group $G=S_n$, the
representations are given by Young diagrams with $n$ boxes. Let
the {\em arm length} $a_k$ of the $k$th box of a Young diagram be the
number of boxes on its right in the same row, and the {\em leg length}
be the number of boxes below it in the same column. The {\em hook
length}
of the $k$th box is then $h_k=a_k+\ell_k+1$. Then the formula for $P(H_m,t)$
can be written as
\[
P(H_m,t)=n!\,t^{mn(n-1)/2}\sum_{\stackrel{\mbox{\tiny{Young}}}
{\mbox{\tiny{diagrams}}}}
\prod_{k=1}^nt^{m(\ell_k-a_k)+\ell_k}\frac
{1-t^k}{h_k(1-t^{h_k})},
\]
see \ref{ss-example} below.

It seems to be impossible to generalize the product formula for $P(H,t)$
for general $m$ since already the dihedral case showed \cite{FeiginVeselov}
that different isotypical components corresponding to
irreducible representations of $G$ behave in a different way as a
function of $m$.

We should mention that the description of the action of the Coxeter
group on the space
of solutions of (\ref{m-har}) can be extracted also from Opdam's
papers \cite{Opdam1},
\cite{Opdam2}\footnote{We are grateful to I. Cherednik
who attracted our attention to these very interesting papers}.
We believe that our derivation is more illuminating and can be used also
for effective description of $m$-harmonic polynomials.

The paper is organised as follows. In Section \ref{s-2}
we review the basic facts about the Calogero--Moser systems and
$m$-harmonic polynomials. The construction of the system of
KZ equations and of the Matsuo--Cherednik isomorphism
is given in Section \ref{s-KZ}. In Section \ref{s-4}
we apply 
this construction to $m$-harmonic polynomials.
We first describe the action of the Coxeter group
on the space of solutions of the KZ equations and give the proof of
Theorem \ref{t-0}. Then we give an explicit 
construction of  $m$-harmonic polynomials
for $G=S_n$ of lowest positive degree using
an integral representation of solutions of the KZ equations.
Finally, in the $S_n$ case,
the asymptotic distribution of degrees of harmonic
polynomials for large $m$ and $n$ are described using results
of Kerov on the statistical properties of large Young
diagrams. We conclude our paper
with Section \ref{s-last}, where we comment on
on interesting recent developments
and open questions in this subject.

\section{Generalized Calogero--Moser systems and $m$-harmonic polynomials}
\label{s-2}

The generalized Calogero--Moser operator related to a Coxeter group $G$
and a multiplicity function $m$ (not necessary integer-valued but
$G$-invariant)
was introduced by Olshanetsky and Perelomov \cite{OP} and has the form
\begin{equation}\label{CaMo}
L=\Delta-\sum_{\alpha\in
A}\frac{m_\alpha(m_\alpha+1)(\alpha,\alpha)}{(\alpha,x)^2}.
\end{equation}
We will use its
gauge equivalent  version ${\mathcal L}= \hat g L \hat g^{-1}$, where
$\hat g$ is the operator of multiplication by
$g=\prod_{\alpha\in A}(\alpha,x)^{m_\alpha}$
which has a form
\begin{equation}\label{CaMog}
\mathcal{L} = \Delta -\sum_{\alpha\in
A}\frac{2m_\alpha}{(\alpha,x)}\partial_\alpha.
\end{equation}

This operator is one of  $n$ commuting operators $\mathcal{L}_1 =
\mathcal{L}, \mathcal{L}_2,...,\mathcal{L}_n$
with the highest symbols $\sigma_i$. One of the best ways to describe
these operators has been discovered
by G. Heckman \cite{Heckman} and uses the following
difference-differential Dunkl operators \cite{Dunkl}
\begin{equation}\label{e-Dunkl}
D_\xi u(x)=\partial_\xi u(x)
+\sum_{\alpha\in A}m_\alpha\frac{(\alpha,\xi)}{(\alpha,x)}(u(s_\alpha x)-u(x)),
\end{equation}
where $\xi \in V$. A remarkable property of these operators is their
commutativity:
$$[D_\xi, D_\eta] = 0$$
for all $\xi, \eta \in V$. This allows us to define the
difference-differential operators
$\sigma_i(D)$ corresponding to the basic invariants $\sigma_i \in S^G$.
Now the operators $\mathcal{L}_i$ can be defined as restrictions of
the operators $\sigma_i(D)$
on the space of $G$-invariant functions (see \cite{Heckman}).

Consider the joint eigenspace of these operators:
\begin{equation}\label{syst}
\begin{cases}
\mathcal{L}_1 \phi= \sigma_1 (\lambda) \phi\\
\hbox to 25mm {\dotfill} \relax\\
\mathcal{L}_n \phi = \sigma_n(\lambda) \phi.
\end{cases}
\end{equation}

\begin{thm}\cite{FeiginVeselov}\label{t-21}
For any Coxeter group $G$ and integer multiplicity
function $m$ all the solutions of the system (\ref{syst}) are
holomorphic everywhere. They form
a space of dimension $|G|$
where the natural action of $G$ is its regular representation.
When $\lambda = 0$ all the solutions are polynomial.
\end{thm}

When all the multiplicities are zero we have
the
usual harmonic polynomials related to a Coxeter group and this result
is well-known (see \cite{Steinberg2,Helgason}).
Following \cite{FeiginVeselov} we will call the solutions of
the system
\begin{equation}\label{syst0}
\begin{cases}
\mathcal{L}_1 \phi= 0\\
\hbox to 15mm {\dotfill} \relax\\
\mathcal{L}_n \phi = 0.
\end{cases}
\end{equation}
{\it $m$-harmonic polynomials} and denote the corresponding
space $H_m$.

The system (\ref{syst}) with generic $m$ has been discussed by
Heckman and Opdam
\cite{HO} (see also \cite{Opdam1})
who showed that it can be represented as a holonomic system
of rank $|G|$ but have never written this system explicitly.
In the next section we present such a representation.\footnote {There is a remark in Opdam's paper
\cite{Opdam1} (see page 335)
which suggests a possibility of such representation.}

\section{Knizhnik--Zamolodchikov equations and the isomorphism theorem}
\label{s-KZ}
In this section, we explain the construction of the Matsuo--Cherednik
isomorphism. It is essentially an adaptation to the rational Coxeter 
group case of a construction of Cherednik (see \cite{Cherednik2}
and \cite{Cheredniklectures}), who treated the trigonometric
case for Weyl groups. However we still present it here, since
in this case the results can be proved more directly, without Hecke
algebra theory. 

We are going first to introduce a certain version of Knizhnik--Zamolodchikov
(KZ) equations on $V$ with values in $G$-modules.
By a $G$-module
we mean a
left module over the group algebra $\C[G]$.
The alternating (sign) representation $\epsilon$ is
defined as the homomorphism $G\to\{1,-1\}$ such that
$\epsilon(s_\alpha)=-1$ for
all reflections $s_\alpha$.
The alternating $G$-module $\C_\epsilon$ is the one-dimensional
module with action given by $\epsilon$.
\subsection{KZ equations with values in $G$-modules}
Let $M$ be a $G$-module, $\pi:V\to \mathrm{End}_\C(M)$ be a $G$-equivariant
linear map 
whose image $\pi(V)$ is contained in a commutative subalgebra. 
The {\em Knizhnik--Zamolodchikov connection} is the system of
first order differential operators
\[
\nabla_\xi=\partial_\xi-\sum_{\alpha\in
A}m_\alpha\frac{(\alpha,\xi)}{(\alpha,x)}
(s_\alpha+1)-\pi(\xi),\qquad \xi\in V,
\]
acting on functions $\psi:V-\cup_\alpha \Pi_\alpha\to M$.

\begin{lemma}\cite{Cherednik, Dunkl2}\label{l-1}
Let $M$ be a $G$-module, $\pi\in\mathrm{Hom}_G(V,\mathrm{End}_\C(M))$, such
that $\pi(\xi)\pi(\eta)=\pi(\eta)\pi(\xi)$ for all $\xi,\eta\in V$. Then
\begin{enumerate}
\item[(i)]
The KZ connection is flat: $\nabla_\xi\nabla_\eta=
\nabla_\eta\nabla_\xi$, $\forall \xi,\eta\in V$.
\item [(ii)] Let $G$ act on $M$-valued functions on $V$ by
$({}^g\psi)(x)=g(\psi(g^{-1}x)), g\in G$.
Then the  KZ connection is $G$-equivariant:
$\nabla_{g\xi}{}^g\psi={}^g(\nabla_\xi\psi)$, $g\in G$.
\end{enumerate}
\end{lemma}

It follows that if $\mathrm{dim}(M)=d$,
the  solutions of the system of
{\em KZ equations with values in $M$},
\[\nabla_\xi\psi(x)=0, \qquad \xi \in V\]
on any connected, simply connected open subset
of $V-\cup_\alpha\Pi_\alpha$ form  a $d$-dimensional vector space.
Moreover, global solutions form a $G$-module.

For $\lambda\in V$, let $I(\lambda)$ be the ideal in $S$ generated
by invariant polynomials vanishing at $\lambda$. It is clearly
invariant under the action of $G$.
We consider the KZ equations with values in the $G$-module
$M(\lambda)=S/I(\lambda)$, with $\pi(\xi)$ being the operator of multiplication
by the linear function $z\mapsto(\xi,z)$:
\begin{equation}\label{e-KZ}
\partial_\xi\psi(x)=
\sum_{\alpha\in A}m_\alpha\frac{(\alpha,\xi)}{(\alpha,x)}
(s_\alpha+1)\psi(x)+\pi(\xi)\psi(x),\qquad \xi\in V.
\end{equation}
We now show that the $G$-module $S/I(\lambda)$
is isomorphic to $\C[G]$. The isomorphism is not canonical and we
choose it in such a way that it is well-behaved as $\lambda\to 0$.
\begin{lemma}\label{l-2}
The composition $H_0\to  S\to M(\lambda)$ of the inclusion with the
canonical projection is an isomorphism of $G$-modules. In particular,
$M(\lambda)$ is isomorphic as a $G$-module to $\C[G]$ with left
action of $G$.
\end{lemma}

\noindent{\it Proof:}
For $\lambda=0$ this is well-known fact about harmonic polynomials
(see \cite{Helgason}).

Let $h_1,\dots, h_{|G|}$ be a basis of $H_0$.
Then, by Chevalley's theorem \cite{Chevalley} (see also \cite{Helgason}),
any polynomial in $S$ may be uniquely written as $q=\sum_{i=1}^{|G|} p_ih_i$,
with $p_i\in S^G$. We claim that $q\in I(\lambda)$ if and only if
$p_i(\lambda)=0$
for all $i$. Indeed, if $p_i(\lambda)=0$ for all $i$ then it is clear
that $q\in I(\lambda)$. Conversely, suppose that
$q=\sum p_i h_i\in I(\lambda)$. Then, by definition, $q=\sum p_k'q_k$ 
with $p_k'\in
S^G$
vanishing at $\lambda$. By expressing each $q_k$ as a linear combination
of $h_i$ with coefficients is $S^G$, get that $q=\sum p_i''h_i$ with
$p_i''\in S^G$ and
$p_i''(\lambda)=0$. By  uniqueness we deduce that $p_i''=p_i$ and the claim
follows.

We thus can define a map of $G$-modules
from $I(0)$ to $I(\lambda)$ by
$\phi_\lambda:\sum p_ih_i\mapsto\sum (p_i-p_i(\lambda))h_i$. It is an
isomorphism with
inverse map $\sum p_ih_i\mapsto \sum (p_i-p_i(0))h_i$.
This map sends polynomials of degree $\leq k$ to polynomials of
degree $\leq k$.
In particular, if we denote by $F^k(S)=\oplus_{j=0}^k S^j$
the polynomials of $S$ of degree
$\leq k$, then $\mathrm{dim}(F^k(S)/(F^k(S)\cap I(\lambda)))$
is independent of $\lambda$
for all $k$, and thus $\mathrm{dim}(S/I(\lambda))=\mathrm{dim}(S/I(0))$.
Moreover, $\phi_\lambda$
maps a homogeneous polynomial $p\in I(0)$ to a polynomial of the form
$p+q$ with $\deg(q)<\deg(p)$.
As a consequence, if $h\neq0$ is in the kernel of
the natural map $H_0\to S/I(\lambda)$, then the component of highest degree of
$h$ is in $I(0)$ and thus vanishes, a contradiction.
\hfill$\square$

\subsection{The 
Matsuo--Cherednik map}
The 
Matsuo--Cherednik map is defined in terms of the function
\[
w_\lambda(y)=\frac
{\sum_{g\in G}\epsilon(g)e^{(g\lambda,y)}}
{\prod_{\alpha\in A}(\alpha,\lambda)}.
\]
\begin{lemma}\label{l-3}
\
\begin{enumerate}
\item[(i)]
The function $w_\lambda(y)$ is holomorphic in $\lambda,y\in V$. It obeys
$w_{\lambda}(y)=\epsilon(g)w_\lambda(gy)=w_{g\lambda}(y)$, $g\in G$.
\item[(ii)]
For any fixed $\lambda$, the function
\[
y\mapsto \frac{w_\lambda(y)}{\prod_{\alpha\in A}(\alpha,y)}
\]
is holomorphic on $V$ and is not zero at $y=0$.
\item[(iii)]
$w_0(y)=C\prod_{\alpha\in A}(\alpha,y)$ for some non zero number $C\in\R$.
\item[(iv)]
Let $\lambda\in V$, $p\in S$. Then $p(\partial)w_\lambda=0$ if and
only if $p\in I(\lambda)$.
\end{enumerate}
\end{lemma}

\noindent{\it Proof:}
(i) Since the numerator obeys $n(s_\alpha\lambda,y)=-n(\lambda,y)$,
$\alpha\in A$,
it is divisible by $\prod_{\alpha\in A}(\alpha,\lambda)$ in the ring
of holomorphic functions.
The behaviour under the $G$-action follows from the fact that the numerator
is skew-invariant as a function of both $\lambda$ and $y$ and the denominator
is a skew-invariant function of $\lambda$. (ii) By the same argument,
$w_\lambda(y)$ is also divisible by $\prod_{\alpha\in A}(\alpha,y)$. Let
$C(\lambda)$ be the value of the quotient at $y=0$. Then the Taylor series
at $0$ of $w_\lambda$ begins with
$w_\lambda(y)=C(\lambda)\prod_{\alpha\in A}(\alpha,y)+\cdots$.
To show that $C(\lambda)\neq 0$ it is sufficient to find a polynomial
$w\in S^N$, $N=|A|$,
such that $w(\partial)w_\lambda(y)$ does not vanish at $y=0$.
Let
$w(z)=\prod_{\alpha\in A}(\alpha,z)$. Then we have
\[
w(\partial)e^{(g\lambda,y)}=w(g\lambda)e^{(g\lambda,y)}=\epsilon(g)
w(\lambda)e^{(g\lambda,y)}.
\]
It follows that $w(\partial)w_\lambda(y)|_{y=0}=|G|\neq0$.
(iii)
It follows from the homogeneity relation
$w_{\lambda}(ay)=a^{N} w_{a\lambda}(y)$,  that $w_0$ is
a homogeneous polynomial of degree $N$. With (ii), this implies the claim.
(iv)
If $p\in S^G$, $p(\partial)w_\lambda=p(\lambda)w_\lambda$.
Thus $p(\partial)w_\lambda$ vanishes if $p\in I(\lambda)$. Conversely,
suppose that $p(\partial)w_\lambda(y)=0$ and let $p_h$ be the homogeneous
term of highest degree in $p$. Then $p_h(\partial)$ vanishes on the
term of lowest degree of $w_\lambda$, namely, by (ii),
$p_h(\partial)\prod_\alpha(\alpha,y)=0$.
It follows, see \cite{Helgason} or Lemma \ref{l-2},
  that $p_h\in I(0)$, i.e.,
$p_h(z)=\sum q_i(z)r_i(z)$
with $r_i\in S^G$ homogeneous of positive degree.
Let $\bar p_h(z)=\sum q_i(z)(r_i(z)-r_i(\lambda))$.
Clearly, $\bar p_h\in I(\lambda)$ and has leading term $p_h(z)$. Thus
$p$ can be
written as
\[
p(z)=\bar p_h(z)+q(z)\equiv q(z) \mod I(\lambda),
\]
where $q$ has lower degree and $q(\partial)w_\lambda=0$.
By induction on the degree, we thus see that $p\in I(\lambda)$.
\hfill $\square$

\medskip

Let $\mu:S\to \C_\epsilon$ be the  $G$-module homomorphism
\[
p\mapsto p(\partial)w_\lambda\left|_{y=0}\right.\, .
\]
By Lemma \ref{l-3}, (iv), $\mu$ induces a map,
also denoted by $\mu$ from $S/I(\lambda)$ to $\C_\epsilon$.

We first formulate the local version of 
the isomorphism theorem.
\begin{thm}\label{t-1}
Let $\lambda\in V$, $U$ a connected, simply connected
open subset of $V-\cap_{\alpha\in A}\Pi_\alpha$. The 
{\em
Matsuo--Cherednik map}
\[
\psi\mapsto \phi= \mu\circ\psi
\]
is an isomorphism from the space
of solutions $\KZ(\lambda;U)$ of the KZ equations \eqref{e-KZ}
on $U$ with values in $S/I(\lambda)$ to the space $\CM(\lambda;U)$
of solutions of the system of differential equations
\begin{equation}\label{e-CM}
\mathcal{L}_i\phi=\sigma_i(\lambda)\phi,\qquad i=1,\dots,n.
\end{equation}
of the Calogero--Moser system (\ref{syst}).

\end{thm}

This theorem is proven in \ref{ss-pt1} below.

\medskip

The global version of the isomorphism theorem follows from
the fact that all local solutions of \eqref{e-CM} extend to global
solutions (see Section \ref{s-2}).

\begin{thm}\label{t-1g}
Let $\lambda\in V$. The  map
\[
\psi\mapsto \phi= \mu\circ\psi
\]
is an isomorphism of $G$-modules from the space
of global solutions $\KZ(\lambda)$ of the KZ equations \eqref{e-KZ}
with values in $S/I(\lambda)$ to the space $\CM(\lambda)\otimes\C_\epsilon$
of global solutions of the system of differential equations
\eqref{e-CM}
of the Calogero--Moser system, with $G$-action twisted
by $\epsilon$.
\end{thm}

Finally, by Theorem \ref{t-21},
all solutions of \eqref{e-CM} extend to holomorphic functions on all
of $V$.  This has the following consequence.

\begin{corollary}\label{c-4}
All solutions of the KZ equations \eqref{e-KZ} extend to holomorphic
functions on $V$. If $\lambda=0$, all solutions of the KZ equations
\eqref{e-KZ} are polynomial.
\end{corollary}

Theorem \ref{t-1g} and Corollary \ref{c-4} are proven in \ref{ss-pt1gc4}.

\medskip

\noindent{\bf Remark.} Matsuo \cite{Matsuo} 
considered 
the KZ equations with values in the left regular $G$-module. He considered
similar maps from the solution space of the  KZ  equations
to $\CM(\lambda)$, which are
isomorphisms as long as $\lambda$ does not belong to the
discriminant locus $\cup_{\alpha\in A}\Pi_\alpha$. 
Cherednik \cite{Cherednik2}, \cite{Cheredniklectures}
 put the construction in the setting of Hecke algebras,
and considered more generally KZ equations with values in representations
of Hecke algebras of Weyl groups (in the trigonometric case). In particular he
studied the KZ equations with values in an induced Hecke algebra
representation and constructed 
a variant of the Matsuo map which is an isomorphism
for all values of the spectral parameters $\lambda$.
In the rational limit, Cherednik's construction reduces to ours (for
Weyl groups).

\subsection{A non-degenerate bilinear form on $S/I(\lambda)$}
A convenient technical tool in the proofs
is a
non-degenerate bilinear form.
Let $O(V)$ be the space of holomorphic functions on $V$.
The bilinear form $\langle p,q\rangle=p(\partial)q(y)|_{y=0}$ on $S$
extends to a  pairing $S\times O(V)\to \C$ defined by the same formula.
For $\lambda\in V$, let $(\ ,\ )_\lambda$ be the symmetric bilinear
form on $S$ defined by
\[
(p,q)_\lambda=p(\partial)q(\partial)w_\lambda(y)|_{y=0}=\langle
p,q(\partial)w_\lambda\rangle.
\]

\begin{lemma} \label{l-ndbf}
The bilinear form $(\ ,\ )_\lambda$ induces a well-defined
symmetric non-degenerate bilinear form on $S/I(\lambda)$.
\end{lemma}

\noindent{\it Proof:}
By Lemma \ref{l-3}, (iv), $(\ ,\ )_\lambda$ vanishes if one of the arguments
is in $I(\lambda)$ so it is well-defined on the quotient. It is clearly
symmetric. Suppose $(p,q)_\lambda=0$ for all $p\in S$. This means
all Taylor coefficients  of $q(\partial)w_\lambda$ vanish. Since it
is a holomorphic function on $V$, this function vanishes identically.
By Lemma \ref{l-3} $q\in I(\lambda)$. This shows that the
induced bilinear form on $S/I(\lambda)$ is non-degenerate.\hfill $\square$

\medskip

Note that the map $\mu$ is $p\mapsto (p,1)_\lambda=\langle p,w_\lambda\rangle$.

\subsection{Proof of Theorem \ref{t-1}}\label{ss-pt1}
We adapt Matsuo's construction to this case. We first show
that the map sends solutions of the KZ equations to eigenfunctions,
following his proof word by word.

We need a description of the polynomial algebra
$\mathcal{D}=\C[\mathcal{L}_1,\dots,\mathcal{L}_n]$
generated by the commuting operators $\mathcal{L}_i$.
These differential operators are defined
through Dunkl operators \cite{Dunkl}, \cite{Heckman},
by the condition that
$\mathcal{L}_i$ acts on invariant functions
as the restriction of the differential-difference
operators $\sigma_i(D)$.
An other system of generators is more convenient here
\cite{Heckman3}.
Define for each $\xi\in V$,
differential operators $D_\xi^{(d)}$ recursively by $D_\xi^{(0)}=1$ and
\[
D^{(d)}_\xi=\partial_\xi D^{(d-1)}_\xi
+\sum_{\alpha\in
A}m_\alpha\frac{(\alpha,\xi)}{(\alpha,x)}(D^{(d-1)}_{s_\alpha\xi}
-D^{(d-1)}_\xi).
\]
Alternatively, $D^{(d)}_\xi$ is the unique differential operator
whose restriction
to $S^G$ coincides with the restriction
to $S^G$ of the $d$th power of the Dunkl differential-difference operator
\eqref{e-Dunkl}.
The algebra $\mathcal{D}$ is generated by the operators
$D_{\xi,d}=\sum_{g\in G}D_{g\xi}^{(d)}$,
$\xi\in V$, $d\in\Z_{\geq0}$. Indeed, $D_{\xi,d}=p_{\xi,d}(D)$
on invariant functions,
where $p_{\xi,d}(x)=\sum_{g\in G}(g\xi,x)$; since any invariant
polynomial can be written as a linear combination of
polynomials $p_{\xi,d}$, it is clear that $\mathcal D$ is
spanned by the operators $D_{\xi,d}$.
Thus the space $\CM(\lambda)$ consists of solutions of the
differential equations
\[
D_{\xi,d}\phi(x)=\sum_{g\in G}(g\xi,\lambda)^d\phi(x)
\]
Let $O(V)$ be the space of holomorphic functions on $V$.

\begin{lemma}\label{l-4.1}
Let $\psi\in\KZ(\lambda;U)$,
$\phi(x)=\mu\circ\psi(x)=\langle\psi(x),w_\lambda\rangle$.
Then
\[
D^{(d)}_\xi\phi(x)=\mu(\pi(\xi)^d\psi(x)).
\]
\end{lemma}

\noindent{\it Proof:} Clearly this holds for $d=0$. Assume inductively that
the claim holds for $D^{(d-1)}_\xi$ and all $\xi\in V$. Then
\begin{eqnarray*}
D^{(d)}_\xi\langle\psi(x),w_\lambda\rangle&=&
\partial_\xi D^{(d-1)}_\xi\langle\psi(x),w_\lambda\rangle
\\ &&+\sum_{\alpha\in
A}m_\alpha\frac{(\alpha,\xi)}{(\alpha,x)}(D^{(d-1)}_{s_\alpha\xi}
-D^{(d-1)}_\xi)\langle\psi(x),w_\lambda\rangle\\
&=&
\partial_\xi\langle\pi(\xi)^{d-1}\psi(x),w_\lambda\rangle
\\ &&+\sum_{\alpha\in A}m_\alpha\frac{(\alpha,\xi)}{(\alpha,x)}
\langle (\pi({s_\alpha\xi})^{d-1}
-\pi(\xi)^{d-1})
\psi(x),w_\lambda\rangle.
\end{eqnarray*}
Since $\psi$ is a solution of the KZ equations, we have
\begin{eqnarray*}
\partial_\xi\langle\pi(\xi)^{d-1}\psi(x),w_\lambda\rangle
&=&
\sum_{\alpha\in A}m_\alpha\frac{(\alpha,\xi)}{(\alpha,x)}
\langle\pi(\xi)^{d-1}(s_\alpha+1)\psi(x),w_\lambda\rangle
\\ &&+\langle\pi(\xi)^{d}\psi(x),w_\lambda\rangle\\
&=&
\sum_{\alpha\in A}m_\alpha\frac{(\alpha,\xi)}{(\alpha,x)}
\langle(-\pi({s_\alpha\xi})^{d-1}+\pi({\xi})^{d-1})\psi(x),w_\lambda\rangle
\\ &&+\langle\pi(\xi)^{d}\psi(x),w_\lambda\rangle,
\end{eqnarray*}
and it follows that
$D^{(d)}_\xi\phi(x)=\langle\pi(\xi)^d\psi(x),w_\lambda\rangle$.
\hfill$\square$

\begin{lemma}
If $\psi\in\KZ(\lambda;U)$ then $\phi(x)=\mu\circ\psi(x)\in \CM(\lambda;U)$.
\end{lemma}

\noindent{\it Proof:} By Lemma \ref{l-4.1},
\begin{eqnarray*}
D_{\xi,d}\phi(x)&=&\langle\sum_{g\in G}(\pi(g\xi))^d\psi(x),w_\lambda\rangle\\
&=&\langle\psi(x),\sum_{g\in G}\partial_{g\xi}^dw_\lambda\rangle.
\end{eqnarray*}
For any invariant polynomial $p\in S^G$, we have
$p(\partial)w_\lambda=p(\lambda)w_\lambda$.
In particular,
\[
\sum_{g\in G}\partial_{g\xi}^dw_\lambda=\sum_{g\in G}(g\xi,\lambda)^dw_\lambda,
\]
which implies the claim. \hfill$\square$

\medskip

It remains to prove that the Matsuo--Cherednik map
$\psi\to\mu\circ\psi$ is an isomorphism
between solution spaces on $U$.
By Theorem \ref{t-21},
$\CM(\lambda;U)$ is a $|G|$-dimensional space.
Since the KZ equations are  a holonomic system
with local solution
space of the same
dimension, it suffices to show that the kernel of the
Matsuo--Cherednik map is trivial.
Let $\psi$ be in the kernel, i.e., $\langle
\psi(x),w_\lambda\rangle=0$, $\forall x\in U$.
By using
the KZ equations, we may express $\langle\psi,\partial_\xi
u\rangle=\langle\pi(\xi)\psi(x),u\rangle$
in terms of $\langle\psi(x),u\rangle$ and $\langle\psi(x),s_\alpha
u\rangle$. It follows
inductively that $\langle\psi(x),p(\partial)w_\lambda\rangle=0$, for
any $p\in S$. Thus,
by Lemma \ref{l-ndbf},
$\psi(x)$ vanishes identically.

\subsection{Proof of Theorem \ref{t-1g} and of Corollary
\ref{c-4}}\label{ss-pt1gc4}
We have
$\mu\circ{}^g\psi(x)=\mu(g\psi(g^{-1}x))=\epsilon(g)\mu\circ\psi(g^{-1}x)$.
This means
that the Matsuo--Cherednik map is a morphism of $G$-modules if we
twist the $G$-action on
$\CM(\lambda)$ by $\epsilon$.
Local solutions of \eqref{e-CM} extend to holomorphic functions on $V$,
see Theorem \ref{t-21}.
Therefore, if $\psi$ is a local solution of the KZ equations, then
$\langle\psi(x),w_\lambda\rangle$ extends to a  holomorphic function
on $V-\cup_\alpha\Pi_\alpha$.
It follows as above, by using the equations, that all coordinates
$\langle\psi(x),p(\partial)w_\lambda
\rangle$ of $\psi$ are meromorphic functions with possible poles on
the reflection hyperplanes.
We now show that these functions are holomorphic on $V$. By Hartogs' theorem it
is sufficient to show that they are holomorphic at generic points of
the hyperplanes.

In the vicinity of a point of the hyperplane $\Pi_\alpha$ which does not lie on
any other hyperplane, a nontrivial solution $\psi$ has a  Laurent expansion
$\psi(x)=(\alpha,x)^k\psi_k
+(\alpha,x)^{k+1}\psi_{k+1}+\cdots$ for some functions
$\psi_k,\psi_{k+1},\dots$  of the
transversal coordinates and $\psi_k\neq0$.
The characteristic exponent $k$ is an eigenvalue of $m_\alpha(s_\alpha+1)$,
and is therefore $0$ or $2m_\alpha\geq 0$. This proves that $\psi$ is regular.

Moreover, we know from Theorem \ref{t-21}
  that if $\lambda=0$,
eigenfunctions $\phi$ are polynomials,
so solutions of the KZ equations are regular rational functions,
i.e.,  polynomials. \hfill$\square$

\section{Applications to $m$-harmonic polynomials}\label{s-4}
In this section, we apply 
the Matsuo--Cherednik isomorphism to the
most degenerate case $\lambda=0$, the case of $m$-harmonic
polynomials.
In particular, we study the action of the Coxeter group on the
space $H_m$ of $m$-harmonic polynomials. As $H_m$ is, as a $G$-module,
isomorphic to the left regular module, we know that the multiplicity
of any simple $G$-module in $H_m$  is  equal to its
dimension. A more detailed information is given by the {\em Poincar\'e
polynomial} of a simple module, which is the generating function
of the multiplicities of the simple module in the space of
$m$-harmonic polynomials of given degree. Let $V_1,\dots,V_p$ be
a list of all inequivalent simple $G$-modules up to isomorphism.
Let
\[
\mu(V_j,H_m^d)=\mathrm{dim}\,\mathrm{Hom}_G(V_j,H^d_m)
\]
be the multiplicity of $V_j$ in the space $H^d_m$ of
homogeneous $m$-harmonic polynomials of degree $d$. We set
\[
P_{j}(H_m,t)=\sum_{d\geq0}\mu(V_j,H_m^d)t^d.
\]
Since $H_m$ is isomorphic to the left regular $G$-module,
the Poincar\'e polynomial of $V_j$ obeys $P_j(H_m,1)=\mathrm{dim}(V_j)$.

Our formula for this polynomial is given in terms of the corresponding
polynomial for $m=0$. Thus we first review what is known for $m=0$.

\subsection{The action of $G$ on harmonic polynomials}
\label{ss-33}
By Chevalley's theorems, the ring of invariants $S^G$ is freely
generated by homogeneous polynomials $\sigma_j$ of degree $d_j$,
$j=1,\dots,n$, and $S$ is a free module over $S^G$. Moreover,
as a basis of $S$ over $S^G$ we may take any basis of the space $H_0$
of harmonic
polynomials.
It follows that the Poincar\'e polynomial $P_j(H_0,t)$  can be expressed
in terms of the
generating series $P_j(S,t)$  of the multiplicities of $V_j$ in the
space of all polynomials of given degree:
\begin{equation}\label{e-PjH}
P_j(H_0,t)=\prod_{k=1}^n(1-t^{d_k})P_j(S,t).
\end{equation}
These generating series  are known for all simple modules
of all Coxeter groups: the case of Weyl groups of
classical root systems is treated by Lusztig in \cite{Lusztig}, \S 2,
(in the $A_\ell$ his formula is deduced from results of
Steinberg \cite{Steinberg}).
Simpler formulae for the classical cases were discovered by Kirillov
\cite{Kirillov}. Formulae for $F_4$, due to Macdonald,  and
tables for the coefficients of $P_j(H_0,t)$ for exceptional
Weyl groups are given in \cite{BeynonLusztig}. As for general Coxeter groups,
the case of the dihedral groups is easy to work out by hand, see the
first example below. The case of $H_4$ was treated in
\cite{AlvisLusztig}. For a comprehensive
overview and a computer program see \cite{GeckPfeiffer}.

\medskip

\noindent{\bf Examples.}\
\begin{enumerate}
\item The dihedral group $G=I_2(N)$ of symmetries of a regular $N$-gon
has exponents $d_1=2$ and $d_2=N$. If $N$ is odd, there are two
one-dimensional modules, the trivial module $V_0$ and the alternating
module $V_N$. There are $N-1$ two-dimensional simple $G$-modules
$V_1=V,V_2,\dots,V_{N-1}$. The $2\pi/N$-rotation has eigenvalues
$\exp(\pm2\pi ij/N)$ in $V_j$.
We then have $P_j(S,t)=P_j(H_0,t)/(1-t^2)(1-t^N)$, and,
if $N$ is odd,
\[
P_0(H_0,t)=1, \quad P_N(H_0,t)=t^N,\quad P_j(H_0,t)=t^j+t^{N-j},
\]
for $1\leq j\leq N-1$.
If $N=2k$ is even, $V_{k}$ is replaced by two one-dimensional modules
$V^{\pm}_k$ with $P_{k}^{\pm}(H_0,t)=t^k$. The other formulae hold without
modification.

\item Let $G=S_n$ be the group of permutations of $n$ letters, generated
by reflections with respect to the hyperplanes $x_j=x_k$ ($j<k$)
of $\R^n$.
Simple $G$-modules
are enumerated by Young diagrams with $n$ boxes. A simple
product formula for $P_j(S,t)$ was given by Kirillov \cite{Kirillov}.
The hook length of a box with coordinates
$(i,j)$ in a Young diagram is the number of boxes
with coordinates $(i,j+p)$ or $(i+p,j)$, $p\geq 0$. Its leg length
is the number of boxes with coordinates $(i+p,j)$, $p>0$. If $h_k$,
$\ell_k$ denote the hook length and the leg length of the $k$th  box of
the Young diagram corresponding to $V_j$, we have
\[
P_j(S,t)=\prod_{k=1}^n\frac{t^{\ell_k}}{(1-t^{h_k})}.
\]
Since the exponents of $S_n$ are $d_j=j$, $j=1,\dots,n$, it follows that
\[
P_j(H_0,t)=\prod_{k=1}^n\frac{t^{\ell_k}(1-t^k)}{(1-t^{h_k})}.
\]
\end{enumerate}

\medskip

\noindent{\bf Remark.}
If $G$ is the Weyl group of a semisimple complex Lie algebra,
$P_j(H_0,t)$ is the generating function of
the multiplicity of $V_j$ in the cohomology groups
of the corresponding flag variety.
Thus  $\sum_j\mathrm{dim}(V_j)P_j(H_0,t)$ is the topological
Poincar\'e polynomial
of the flag variety.

\subsection{Coxeter groups and Hecke algebras} \label{subs-hecke}
The calculation of the action of $G$ on the space of solutions of
the Knizhnik--Zamolodchikov equations uses the fact (discovered in
a special case in \cite{TsuchiyaKanie}, and generalized in \cite{Cherednik1},
see also \cite{Heckman2} for the result in the
context of Calogero--Moser systems) that
the monodromy representation of this system of differential
equation factors through a Hecke algebra.

  Let $(u_\alpha)_{\alpha\in A}$ be
indeterminates associated to reflections of a Coxeter group $G$ so that
$u_\alpha=u_\beta$ if the corresponding reflections are conjugated in $G$.
Fix a presentation of $G$ with set of generators $\Sigma$ and relations
$s^2=1$, $(st)^{m(s,t)}=1$, $s,t\in \Sigma$. We write $u_s=u_\alpha$
if $s=s_\alpha\in\Sigma$.
Accordingly, we write $m_s=m_\alpha$ if
$s=s_\alpha\in\Sigma$.
The length $\ell(g)$ of $g$ is the minimal
number of factors needed to write $g$ as a product of generators.
The Iwahori--Hecke algebra
$H(G)$ of the Coxeter group $G$ is defined as the algebra over
$R=\Z[u_s,s\in\Sigma]$ spanned by elements $T_g$,
$g\in G$ and multiplication rules $T_sT_g=T_{sg}$ if $\ell(sg)>\ell(g)$
and $T_sT_g=u_s T_{sg}+(u_{s}-1)T_g$  if $\ell(sg)<\ell(g)$, $s\in \Sigma$.

If we specialise the parameters
$u_s$ to $1$ we recover the group algebra of $G$. The algebra $H(G)^{K}
=H(G)\otimes_R K$ over the field $K=\C(\sqrt{u_s},s\in\Sigma)$ of
rational functions of  $\sqrt{u_s}$, $s\in\Sigma$,
is split semisimple and isomorphic to the group algebra $KG$.
Moreover every irreducible character $\chi_j$
of $G$ is the specialization
at $\sqrt{u_s}=1$ of an irreducible character
$\hat\chi_j(\sqrt{u_s}\,)$ of $H(G)^K$
with values in $\C(\sqrt{u_s}\,)$
(Benson--Curtis,
Lusztig, Alvis--Lusztig, Digne--Michel,
see \cite{GeckPfeiffer}, Theorem 9.3.5).

The function $\mu\mapsto\hat\chi_j(\exp(\pi i\mu_\alpha))$ defined
for $G$-invariant collections $\mu=(\mu_\alpha)_{\alpha\in A}$ of
complex numbers is a holomorphic function whose values at integer points
is the character of an irreducible $G$-module and we
have for $\mu_s=m_s\in \Z$,
\[
\hat\chi_j(\exp(\pi im_s))=\chi_{\pi_m(j)},
\]
for some permutation $\pi_m$ of the set of equivalence classes of
irreducible $G$-modules. Clearly $\pi_{m+m'}=\pi_m\circ\pi_{m'}$,
$\pi_0=\mathrm{id}$ and $\pi_m$ depends only on the class of the integers
$m_\alpha$ in $\Z/2\Z$.
In particular it is an involution.
We will need two simple properties of the homomorphism $\pi$:
\begin{enumerate}
\item[(A)] Let $j\mapsto j^*$ be the tensor multiplication by the
alternating representation. Then $\pi_m(j)^*=\pi_m(j^*)$.
\item[(B)] $\chi_{\pi_m(j)}(s)=\chi_{j}(s)$ for reflections $s$.
\end{enumerate}
The first property follows from the fact that the automorphism
$\tau$ of $H(G)$
such that $\tau(T_s)= -T_s+(1-u_s)1$ for reflections $s$,
induces an involution $j\mapsto j^*$
on the
set of irreducible modules over $H(G)^K$, which specialises to
the tensor multiplication by the alternating representation
at $u_s=1$.
The character $\chi_{{j^*}}$ is then the specialisation at
$\sqrt{u_s}=1$ of a character $\hat\chi_{j^*}$
with values in $\C(\sqrt{u_s})$.
By using the fact that, for any $g\in G$,
$\tau(T_g)\in H(G)$ has coefficients
which are polynomials in $u_s$, we have
\[
\chi_{\pi_m(j)}^*(g)=\hat\chi_{j}(\sqrt{u_s}\,)
(\tau(T_g))\left|_{\sqrt{u_s}=e^{\pi im_s}}\right.
=\hat\chi_{j^*}(\sqrt{u_s}\,)
(T_g)\left|_{\sqrt{u_s}=e^{\pi im_s}}\right.
=\chi_{\pi_m(j^*)}(g).
\]
This proves the first property.

The second property follows from the relation $(T_s+1)(T_s-u_s)=0$
for reflections $s$
(take $g=s$ in the multiplication rule). This relation
implies that the character of any simple module evaluated at
$T_s$ has the form $a_s u_s-b_s$, where $a_s$, $b_s$ are the multiplicities
of the eigenvalues $u_s$, $-1$. So there are no square roots here and
$\pi_m$ acts trivially.

The explicit description of the
homomorphism $m\mapsto \pi_m$
can be extracted from the literature on Hecke algebras and was
given in \cite{Opdam2}. In fact, $\pi_m$ is trivial for $A_n,B_n,D_n,E_6,
F_4$ and the dihedral groups of order $2N$, with $N$ odd. For dihedral
groups of order $4n$, $n\geq 3$, (this includes $G_2$),
$\pi_m$ is trivial on one-dimensional
modules; on two-dimensional modules it
is the tensor multiplication by the one-dimensional representation $\chi_m$
such that $\chi_m(s_\alpha)=(-1)^{m_\alpha}$. In the remaining cases
$E_7,E_8,H_3,H_4$, for which the reflections are all in the
same conjugacy class, $\pi_{m=1}$ exchanges one or two pairs of
modules (see \cite{Opdam2} for details)
and is the identity on all others.

\subsection{The action of $G$ on solutions of the KZ equations
and on $m$-harmonic polynomials}
Let $C_1,\dots, C_q$ be the conjugacy classes of reflections in $G$.
The multiplicity function $m$ may be regarded as a function on the
set of reflections which is constant on each conjugacy
class.
We denote accordingly by $m_a$ the value of $m$ on $C_a$. We also
let $N_a$ be the number of elements of $C_a$.
\begin{thm}\label{t-5}  Let $V_{j,\alpha}^\pm=
\{v\in V_j\,|\, s_\alpha v=\pm v\}$
be the $(\pm1)$-eigenspace of $s_\alpha$ acting on the simple $G$-module
$V_j$ and set
\begin{equation}\label{e-defd2}
d^\pm_a(V_j)=\frac{2N_a\mathrm{dim}(V_{j,\alpha}^\pm)}{\mathrm{dim}(V_j)},\qquad
a=1,\dots,q,
\end{equation}
for any $s_\alpha\in C_a$. Then
\begin{equation}\label{e-Carlosonioetamo}
P_{\pi_m(j)}(H_m,t)=t^{\sum_{a=1}^qm_ad^-_a(V_j)}P_{j}(H_0,t).
\end{equation}
\begin{equation}\label{e-Cortigianirazzadannata}
P_{j}(H_m,t)=t^{\sum_{a=1}^qm_ad^-_a(V_j)}P_{\pi_m(j)}(H_0,t).
\end{equation}
\end{thm}

Note that \eqref{e-Cortigianirazzadannata} is an immediate
consequence of \eqref{e-Carlosonioetamo}, the fact that
$\pi_m$ is an involution and property (B) of $\pi_m$.

This Theorem can also be obtained as a
special case of a result of Opdam
\cite{Opdam1, Opdam2}, see in particular
Theorem 7 in \cite{Opdam2}.
We give a more transparent and effective proof
in \ref{ss-43} below, based on the computation of the action of
$G$ on the space of solutions of the KZ equations.

{}From the formula \eqref{e-Carlosonioetamo}
we obtain a formula for the total Poincar\'e polynomial
$P(H_m,t)=\sum_{d\geq0} \dim(H_m^d)t^d$.
It is obtained by summing the
contribution of the individual $G$-modules, by
  noticing that $\pi_m$
is a permutation preserving the dimension of the modules:
\[
P(H_m,t)=\sum_j\mathrm{dim}(V_j)\,t^{\sum_{a=1}^qm_ad^-_a(V_j)}P_{j}(H_0,t).
\]

\begin{corollary} We have Poincar\'e duality:
\[
P_{j^*}(H_m,t)=t^{M}P_j(H_m,t^{-1})
\]
The total Poincar\'e polynomial $P(H_m,t)$ has degree
$M=\sum_{\alpha\in A}(2m_\alpha+1)$ and obeys
\[ P(H_m,t)=t^MP(H_m,t^{-1})\]
\end{corollary}

\noindent{\it Proof:} The first formula
follows from Poincar\'e duality for $H_0$, see \eqref{poinc},
the fact that $\pi_m(j^*)=\pi_m(j)^*$ (Property (A)), and
the identity $d_a^-(V_{j^*})=2N_a-d_a^-(V_j)$. The
formula for the total Poincar\'e polynomial is
then an easy consequence. From this formula, we deduce
that $P(H_m,t)$ is a polynomial of degree $\leq M$. The
alternating module (and only this module)
gives a term  $t^M$. Thus the polynomial is of degree
$M$.
$\square$

The proof of Theorem \ref{t-0} is complete.

We now show that $d^-_a$ may also be calculated in terms of $P_j(H_0,t)$,
generalizing a result of Solomon \cite{Solomon}.
The group $\mathrm{Hom}(G,\C^\times)$
of one-dimensional representations of $G$ is
generated by representations $\chi_a:G\to\{1,-1\}$, $a=1,\dots,q$ such that
\[
\chi_a(s)=\left\{\begin{array}{rl}
-1,& \text{if $s\in C_a$,}\\
1,& \text{if $s\in C_b$, $b\neq a$,}
\end{array}
\right.
\]
The group $\mathrm{Hom}(G,\C^\times)$ acts on the set of simple $G$-modules
by tensor multiplication. In particular, $\chi_a$ acting on $V_j$ gives
a module isomorphic to $V_{j\otimes a}$ for some permutation
$j\mapsto j\otimes a$
of $\{1,\dots,p\}$.

\begin{proposition}\label{t-6}
\[
d^-_a(V_j)=N_a+
\left.\frac{d}{dt}\right|_{t=1}
\ln\frac{P_j(H_0,t)}{P_{j\otimes a}(H_0,t)},\qquad
a=1,\dots, q.
\]
If $q=1$ the formula simplifies to
\[
d^-(V_j)=d^-_1(V_j)=2\left.\frac{d}{dt}\right|_{t=1}\ln{P_j(H_0,t)}.
\]
\end{proposition}

\noindent{\bf Remark.} The second formula in Theorem \ref{t-5} is due
to Solomon \cite{Solomon}.
It can be stated, following his formulation, in the following
equivalent  way: if $V_j$ appears in the decomposition of $H_0^q$
into simple $G$-modules for
$q=q_1,\dots,q_{\mathrm{dim}(V_j)}$ then $d^-(V_j)/2$
is the average of the $q_i$. From our result follows the non-obvious
fact that this average is always half an integer.

\subsection{Action of $G$ on $\KZ(0)$ and proof of Theorem
\ref{t-5}}\label{ss-43}
The proof is based on the $G$-module isomorphism between the space of
$m$-harmonic
polynomials and the space $\KZ=\KZ(0)$ of solutions of the KZ equations with
values in $M=M(0)=S/I(0)$.
The first observation is that the space $M$ of values has a filtration
preserved by the KZ connection:
let $F_d(S)=\oplus_{j\geq d}S^j$ be the subspace of polynomials of
degree $\geq d$ and set
$F_d(M)=F_d(S)/(F_d(S)\cap I(0))$. Since $I(0)$ is the direct sum of its
homogeneous components, and the action of $G$ preserves degrees,
we have natural inclusion maps of
$G$-modules,
\[
\C w_0(z)=F_N(M)\subset F_{N-1}(M)\subset\cdots\subset F_0(M)=M.
\]
Here, as usual, $N=\deg w_0(z)$ is the number of reflections of $G$.
The KZ connection preserves functions with values in each
$F_j(M)$, and we thus have a corresponding $G$-module filtration
of the space $\KZ$:
\[
\C w_0(z)=\KZ_N\subset \KZ_{N-1}\subset\cdots\subset\KZ_0=\KZ.
\]
The space $\KZ_N$ is spanned by the constant solution $x\to w_0(z)\in M$,
which is mapped to a constant eigenfunction under the Matsuo--Cherednik map.

Additionally to this filtration, we also have a grading with respect to
the total degree: if we represent
a polynomial function $x\to\psi(x)$ on $V$ with values in $M$ by a scalar
function $\psi(x,z)$ of two variables,
we say that $\psi$ has {\em total degree} $\delta$ if
$\psi(\lambda x,\lambda^{-1}z)=\lambda^\delta\psi(x,z)$. Then the KZ connection
decreases the total degree by $1$. Therefore we may decompose the space of
solutions in homogeneous components with respect to the total degree:
\[
\KZ_{d}=\oplus_{\delta\geq-N}\KZ_{d}^\delta.
\]
\begin{lemma}\label{l-7}
The Matsuo--Cherednik map sends solutions of total degree $\delta$ to
homogeneous $m$-harmonic polynomials of degree $\delta+N$.
\end{lemma}

\noindent{\it Proof:}
The Matsuo--Cherednik map sends $\psi(x,z)$ to
$\psi(x,\partial_y)w_0(y)|_{y=0}$,
which is the value at $y=0$ of a homogeneous polynomial in $x,y$ of degree
$\delta+N$. \hfill $\square$

\medskip

We are thus left to describe solutions in $\KZ_d^\delta$. Suppose
$\psi\in KZ_d$
has a nontrivial projection in ${\KZ}^\delta_d/\KZ^\delta_{d+1}$.
Then $\psi=\psi^d+\psi^{d+1}+\cdots+\psi^N$ with
$\psi^j$ taking values in $M^j=S^j/(S^j\cap I(0))$ . Then the {\em lowest
component} $\psi^d$ is a nontrivial solution of the KZ equations
with values in $M^d$ without $\pi$-term:
\begin{equation}\label{e-KZ0}
\partial_\xi\psi^d(x)=
\sum_{\alpha\in A}m_\alpha\frac{(\alpha,\xi)}{(\alpha,x)}
(s_\alpha+1)\psi^d(x),\qquad \xi\in V.
\end{equation}
\begin{lemma} The map $\psi\to\psi^d$ is an isomorphism of
$G$-modules from $\KZ_d^\delta/\KZ_{d+1}^\delta$ onto the
space of homogeneous polynomial
solutions of \eqref{e-KZ0} of degree $\delta+d$
with values in $M^d=S^d/(S^d\cap I(0))$,
homogeneous of degree $\delta+d$.
\end{lemma}

\noindent{\it Proof:}
The map $\psi\mapsto\psi^d$ is an injective $G$-module homomorphism
from $\KZ_d/\KZ_{d+1}$ to the space of polynomial solutions of \eqref{e-KZ0}
with values in $M^d$. Taking the direct sum over $d$ and comparing dimensions
we see that all solutions of \eqref{e-KZ0} are polynomial and that the
map is an isomorphism.
It maps
$\KZ^\delta_d/\KZ^\delta_{d+1}$ onto the space of solutions $\psi$ of
\eqref{e-KZ0} with $\psi(\lambda x,z)=\lambda^d\psi(\lambda x,\lambda^{-1}z)=
\lambda^{\delta+d}\psi(x,z)$.
\hfill$\square$

\medskip

To describe solutions of \eqref{e-KZ0} with values in $M^d$ we fix
a decomposition of $M^d$ as a direct sum of simple modules and notice
that both sides of \eqref{e-KZ0} contain operators preserving
this decomposition.

\begin{proposition}
Let $\KZ(V_j)$ be the space of  polynomial solutions of
\eqref{e-KZ0} with values
in a simple $G$-module $V_j$. Then
\begin{enumerate}
\item[(i)] Solutions in $\KZ(V_j)$ are homogeneous polynomials
of degree $\sum_{a=1}^qm_ad^+_a(V_j)$, see \eqref{e-defd2}.
\item[(ii)] $\KZ(V_j)$ with action $\psi\mapsto {}^g\psi$,
${}^g\psi(x)=g\cdot\psi(g^{-1}x)$, $g\in G$, is isomorphic
to $V_{\pi_m(j)}$ as a $G$-module (see \ref{subs-hecke}).
\end{enumerate}
\end{proposition}

\noindent{\it Proof:}
(i) Let $E\psi(x)=\frac d{dt}\psi(tx)|_{t=1}=\sum_{j=1}^n
x_j\partial\psi(x)/\partial x_j$
be the Euler vector field. The KZ equation implies
\[
E\psi(x)=\sum_{\alpha\in A}m_\alpha(s_\alpha+1)\psi(x).
\]
The operator on the right-hand side is a central element of $\C[G]$
(it clearly commutes with all reflections) and is thus, by Schur's lemma,
a scalar $a$ times the identity. To compute $a$ we take the trace
and obtain
\begin{eqnarray*}
a&=&\frac1{\mathrm{dim}(V_j)}\sum_{\alpha\in A}m_\alpha \mathrm{tr}_{V_j}
(s_\alpha+1)
\\
&=&
\frac1{\mathrm{dim}(V_j)}\sum_{\alpha\in A}m_\alpha 2\,\mathrm{dim}\,
\mathrm{Ker}(s_\alpha-1:V_j\to V_j).
\end{eqnarray*}
Since the terms of the sum with $s_\alpha$ in the same conjugacy class
are equal to each other, we obtain the claimed formula.

(ii) We apply the method of \cite{Opdam2}, which consists of considering
the KZ equation with complex $G$-invariant
multiplicities $\mu_\alpha$. Let $\psi_\mu$ be
the  solution of the KZ equation
\[
\partial_\xi\psi(x)=
\sum_{\alpha\in A}\mu_\alpha\frac{(\alpha,\xi)}{(\alpha,x)}
(s_\alpha+1)\psi(x),\qquad \xi\in V.
\]
on some neighborhood of some regular point $x_0$
with values in $\mathrm{End}(V_j)$ and
initial condition $\psi_\mu(x_0)=\mathrm{id}$.
If $g\in G$ and
$\gamma$ is a path in $V-\cup_\alpha\Pi_\alpha$ from $x_0$ to $g^{-1}x_0$, we
set
${}^g\psi_\mu(x)=g\psi_\mu(g^{-1}x)$, where the value at $g^{-1}$ is obtained
by analytic continuation along $\gamma$. Then ${}^g\psi_\mu=\psi_\mu\circ T_g$
for some monodromy matrix $T_g$. By the results of  \cite{Cherednik1}, the
$T_g$ obey, for generic $\mu$,
the relations of the Iwahori--Hecke algebra with $u_\alpha=\exp(2\pi
i \mu_\alpha)$.
At $\mu=0$, $T_g$ is given by the action of the group on $V_j$. More generally,
for all integer $\mu$,
$g\mapsto T_g$ is a representation of $G$, since we know that all solutions are
polynomial. It follows that the character of the Iwahori--Hecke
algebra corresponding
to the representation $g\mapsto T_g$ at generic $\mu$ coincides with the
character
$\hat\chi_j(\exp(\pi i\mu_\alpha,\alpha\in A))$ described in 
\ref{subs-hecke} and
the claim follows.
\hfill$\square$

\medskip

We are ready to complete the proof. Solutions of \eqref{e-KZ0} with
values in a simple module $V_j$ occurring in the decomposition of $M^d$
give rise to a subspace of solutions in $\KZ^\delta_d$ of total degree
$\delta=\sum_am_ad^{+}_a(V_j)-d$ and isomorphic to $V_{\pi_m(j)}$. 
The Matsuo--Cherednik
map sends these solutions to $m$-harmonic polynomials of degree $\delta+N$
that form a $G$-module isomorphic to 
$V_{\pi_m(j)^*}=V_{\pi_m(j)}\otimes\C_\epsilon$.
On the other hand, $M^d$ is isomorphic to the space of harmonic polynomials
of degree $d$. Thus
\begin{eqnarray*}
P_{\pi_m(j)^*}(H_m,t)&=&t^{\sum_{a=1}^q m_ad^+_a(V_{j})+N}P_j(H_0,t^{-1})
\\
&=&t^{\sum_{a=1}^q m_ad^-_a(V_{j^*})}P_{j^*}(H_0,t),
\end{eqnarray*}
by Poincar\'e duality
for $H_0$
and the fact that $d^+_a(V_{j})=d^-_a(V_{j^*})$.

The proof of Theorem \ref{t-5} is complete.

\subsection{Example: Lowest degree $m$-harmonic polynomials for $G=S_n$}
\label{ss-example}
If $G=S_n$, Kirillov's formula gives
\[
d^-(V_j)=\sum_{k=1}^n(\ell_k-a_k)+n(n-1)/2
\]
where $a_k=h_k-\ell_k-1$ is the ``arm length'' of the $k$th box of the
Young diagram corresponding to $V_j$. Therefore,
\[
P_j(H_m,t)=t^{mn(n-1)/2}\prod_{k=1}^nt^{m(\ell_k-a_k)+\ell_k}\frac
{1-t^k}{1-t^{h_k}}.
\]
Thus, a part from the trivial module, which occurs in degree 0,
the lowest degree $m$-harmonic polynomials appear in degree $nm+1$.
They form an $(n-1)$-dimensional simple $G$-module
corresponding to the minimal leg-length, maximal arm-length
partition $(n-1,1,0,\dots,0)$. It is isomorphic to the
module $\{x\in\C^n\,|\,\sum x_i=0\}$ with permutation action of $S_n$.

These $m$-harmonic polynomials are associated to solutions of the
KZ equations with values in $F_{N-1}(M)=M^N\oplus M^{N-1}$,
$M=S/I(0)$, $N=n(n-1)/2$.
Let $\psi$ be such a solution,
$\phi(x)=\langle\psi(x),w_0\rangle$ the corresponding
$m$-harmonic polynomial, and $\psi_k(x)=\langle\psi(x),\partial_k w_0\rangle$,
where $\partial_k$ is the derivative with respect to the $k$th
coordinate of $V=\C^n$. Then $\psi$ is uniquely determined by the
components $\phi,\psi_1,\dots,\psi_n$. They obey
\begin{equation}\label{e-sum}
\psi_1+\cdots+\psi_n=0.
\end{equation}
In these terms, the KZ equations are equivalent to the system
\begin{eqnarray}
\label{e-phi}\frac{\partial\phi}{\partial x_i}&=&\psi_i,\\
\label{e-psi}\frac{\partial\psi_i}{\partial
x_j}&=&-m\frac{\psi_i-\psi_j}{x_i-x_j},
\qquad i\neq j,
\end{eqnarray}
(the equations for $\partial_i\psi_i$ are redundant in view of \eqref{e-sum}).

The simplest solutions of these equations are the ones with $\psi_i=0$.
They give rise to constant $m$-harmonic polynomials $\phi=\mathrm{const}$,
which form a trivial $G$-module.

We next describe the solutions corresponding to $m$-harmonic polynomials
of degree $mn+1$. We assume that $m\geq 1$ (the case $m=0$ is left  as
an exercise).
It is easy to check that for any family of
cycles $\gamma(x)$ in the relative integral homology
$H_1(\C,\{x_1,\dots,x_n\})$,
the following formula gives a solution of the
equations \eqref{e-psi}
\[
\psi_j(x)=\int_{\gamma(x)}
\frac{dt}{t-x_j}\prod_{k=1}^n(t-x_k)^{m}.
\]
Indeed the integrand obeys the equations, and the derivative with respect
to the endpoints appearing in $\partial_i\psi_j$ with
$i\neq j$ vanishes since the integrand vanishes there for $m\geq1$.
Moreover, we have
\[
\sum_{j=1}^n\psi_j(x)=\frac1m\int_{\gamma(x)}d \prod_{k=1}^n(t-x_k)^m=0,
\]
since $m\geq 1$.
The corresponding $m$-harmonic polynomials of degree $mn+1$
are obtained by integrating
\eqref{e-phi}:
\[
\phi(x)=\frac1{1+nm}\sum_{i=1}^nx_i\psi_i(x)=\frac1{1+nm}\int_{\gamma(x)}
\sum_{i=1}^n\frac{x_idt}{t-x_i}\prod_{k=1}^n(t-x_k)^{m}.
\]
It remains
to show that we get $n-1$ linearly independent solutions in this way.
A basis of $H_j(\C,\{x_1,\dots,x_k\})$ for generic $x$ is given by
paths  $\gamma_j$ joining $x_1$ to $x_j$, $2\leq j\leq n$. Let
$\phi^{(j)}$ be the corresponding $m$-harmonic polynomials:
\[
\phi^{(j)}(x)=\frac1{1+nm}\int_{x_1}^{x_j}
\sum_{i=1}^n\frac{x_idt}{t-x_i}\prod_{k=1}^n(t-x_k)^{m},
\qquad j=2,\dots,n
\]
These functions are linearly independent, as can be seen by
comparing their behaviour at infinity:
\[
\lim_{x_j\to\infty}x_j^{-mn-1}\phi^{(k)}(x)=(-1)^{m-1}\frac{(m-1)!^2}{
(1+nm)(2m-1)!}
\delta_{jk}.
\]

\subsection{Proof of Proposition \ref{t-6}}

Let $g^k$ denote the action of $g\in G\subset GL(V)$ on the symmetric
power $S^k$ of
the reflection module $V$. Then the character of $S^k$ is $\chi_{S^k}(g)=
\mathrm{tr}(g^k)$
and $\sum_{k=0}^\infty\mathrm{tr}(g^k)t^k=\mathrm{det}^{-1}(1-gt)$ is
the generating
series of these characters. If $\chi_j$ is the character of $V_j$,
the  orthogonality relations of characters implies that
\begin{eqnarray*}
P_j(S,t)&=&\frac1{|G|}\sum_{g\in G}\chi_j(g^{-1})\frac1{\mathrm{det}(1-gt)}\\
         &=&\frac1{|G|}\sum_{g\in
G}\chi_j(g^{-1})\prod_{i=1}^n\frac1{1-\lambda_i(g)t}.
\end{eqnarray*}
Here $\lambda_1(g),\dots,\lambda_n(g)$ are the eigenvalues of $g\in GL(V)$.
It follows from \eqref{e-PjH} that
\begin{eqnarray*}
P_j(H_0,t)&=&
\frac1{|G|}
\sum_{g\in G}\chi_j(g^{-1})\prod_{i=1}^n\frac{1-t^{d_i}}{1-\lambda_i(g)t}\\
&=& \frac1{|G|}
R(t)\prod_{i=1}^n\frac{1-t^{d_i}}{1-t},\\
R(t)&=&
\sum_{g\in G}\chi_j(g^{-1})\prod_{i=1}^n\frac{1-t}{1-\lambda_i(g)t}.
\end{eqnarray*}
If we take  the  derivative of $R(t)$ at $t=1$, the only terms
surviving in the sum
over $G$ are  the
reflections, since they are the only group elements with multiplicity
of the eigenvalue $1$ equal to $n-1$. On the other hand,
$R(1)=\mathrm{dim}(V_j)$
since only the identity survives in this case.
We then have
\begin{eqnarray*}
\left.\frac{d}{dt}\right|_{t=1}\ln P_j(H_0,t)&=&
\sum_{i=1}^n\left.\frac{d}{dt}\right|_{t=1}\ln\frac{1-t^{d_i}}{1-t}
+\frac1{\mathrm{dim}(V_j)}
\sum_{\alpha\in
A}\left.\frac{d}{dt}\right|_{t=1}\frac{1-t}{1+t}\chi_j(s_\alpha)\\
&=&
\sum_{i=1}^n\frac{d_i-1}2+
\frac1{\mathrm{dim}(V_j)}
\left(-\frac12\right)\sum_{\alpha\in A}\chi_j(s_\alpha)
\end{eqnarray*}
Let $\chi_a$ be the one-dimensional representation of $G$
so that $\chi_a(s)=-1$ only on reflections $s\in C_a$.
Since $\chi_{j\otimes a}(s_\alpha)=\chi_a(s_\alpha)\chi_j(s_\alpha)$,
we obtain
\[
\left.\frac{d}{dt}\right|_{t=1}\ln \frac{P_j(H_0,t)}{P_{j\otimes a}(H_0,t)}=
-\frac1{\mathrm{dim}(V_j)}
\sum_{s_\alpha\in C_a}\chi_j(s_\alpha).
\]
Since $\chi_j(s_\alpha)=\mathrm{dim}(V_j)-2\,\mathrm{dim}(V^-_{j,\alpha})$,
the result follows.
The second formula in Theorem \ref{t-6}, valid for groups
with one conjugacy class of reflections, is obtained from the first
formula by applying Poincar\'e duality $P_{j\otimes\epsilon}(H_0,t)=
t^{N}P_{j}(H_0,t^{-1})$. In this case the only one-dimensional
representation is the sign $\epsilon$.

\subsection{Distribution of the degrees of $m$-harmonic polynomials 
for large $m$}

Let us consider again the case of the symmetric group $G = S_n$.
Choose any basis in the corresponding space of $m$-harmonic polynomials $H_m$.
Let $d_1^{(m,n)}=0, d_2^{(m,n)}=d_3^{(m,n)}=...=d_n^{(m,n)}=mn+1, ...,\,
d_{n!}^{(m,n)}= (2m+1)n(n-1)/2$ be the degrees of 
the $m$-harmonic polynomials of the basis
written in increasing order (with repetitions).
We are interested in the distribution of these degrees for large $m$.

As we have shown above all the degrees $d_j^{(m,n)}$ are growing linearly
as a function of $m$ so it is natural to consider the limits of the 
scaled degrees
$\bar d_j^n = \lim_{m\to\infty}d_j^{(m,n)}/m$.
 From our Theorem \ref{t-0} it follows that the set of $\bar d_j^n$
coincide with the set of the numbers $d^{-}(U), U \in \hat S_n$ where 
$\hat S_n$
is the set of all irreducible representations of $S_n$. Moreover each 
number $d^{-}(U)$
appears $(\dim U)^2$ times since the action of $G$ on $H_m$ is 
regular and therefore
each representation $U$ appears with multiplicity equal to its dimension.

On the set $\hat S_n$, which can be represented as
the set of all Young diagrams $Y_n$ of size $n$, there exists a natural
probability measure called {\it Plancherel measure}: each 
representation $U$ has a weight
equal to $(\dim U)^2/n!$ (see e.g. \cite{VeK1}). The function $d^{-}(U)$
can be considered as a random variable on this set.
Let us consider the degree $d_j^{(m,n)}$ of $m$-harmonic polynomial
as a random variable on the discrete set of
$j \in {1,2,...,n!}$ where each $j$ has the same probabilitiy $1/n!$ 
(in other words
we assume that each of the basic $m$-harmonic polynomials has the same weight).
We have the following

\begin{proposition}\label{4.7}
The distribution of the scaled degrees $d_j^{(m,n)}/m$ of 
$m$-harmonic polynomials
as $m$ goes to infinity converges to the distribution of the random 
variable $d^{-}$ on the set
of Young diagrams of size $n$ considered with Plancherel measure.
\end{proposition}

In 1977 Vershik and Kerov \cite{VeK1} and independently Logan and 
Shepp \cite{LoS}
have found remarkable results on the random Young diagrams of size $n$
for large $n$. They have shown that as $n$ goes to infinity the random shape of
such a diagram after proper rescaling concentrates near some universal curve.
The Gaussian fluctuations around this limiting shape have been described
later by Kerov \cite{Ke1} (a
detailed proof of Kerov's result can be found in
recent paper by V. Ivanov and G. Olshanski  \cite{IO}).

We can use 
Kerov's following result to describe the distribution of
the random variable $d^{-}(U)$ on $Y_n$ for a large $n$.
Let us introduce following \cite{IO} the normalised characters
$$p_2^{(n)}(U) =\frac{n(n-1) \chi_U(s)}{\dim U} = \frac {n(n-1) (\dim 
U^+ - \dim U^-)}{\dim U},$$
where $\chi_U$ is the character of $U$, $s = s_{ij} \in S_n$ is any 
transposition and
$U^+, U^-$ are the $\pm$-eigenspaces of its action on $U$.

\medskip

\noindent{\bf Kerov's theorem \cite{IO,Ke1}.}
{\it When $n$ goes to infinity
the distribution of the random variables $\frac{p_2^{(n)}}{\sqrt 2 
n}$ converges to the standard normal distribution
$N(0,1)$.}

\medskip
Since $p_2^{(n)}(U) =  \frac {n(n-1) (\dim U^+ - \dim U^-)}{\dim U} =
n(n-1) - 2d^{-}(U)$
we have as a corollary that
when $n$ goes to infinity the distribution of
$\frac{1}{n}(d^{-} - \frac{n(n-1)}{2})$ converges to $N(0,1/2),$
where $N(0,1/2)$ is the normal distribution
with mean 0 and variance 1/2.

This gives some idea of how the degrees of $m$-harmonic polynomials 
of the symmetric group $S_n$
are distributed for large $m$ and $n$.
We should note that Proposition \ref{4.7} has an obvious generalization for 
any Coxeter group
but Kerov's result is known only for symmetric groups.

\section{Concluding remarks}\label{s-last}

The space $H_m$ of $m$-harmonic polynomials has been introduced in 
\cite{FeiginVeselov} in relation
with the algebra of quasiinvariants $Q_m$ of the Coxeter group $G$. 
In \cite{FeiginVeselov} it was shown that
if $G$ is any dihedral group and the multiplicity $m$ is a constant 
function then
$Q_m$ is freely generated by any basis in $H_m$ as a module over the 
subring $S^G$ of $G$-invariant
polynomials and it was conjectured that the same is true in general situation.
However this turned out not to be true as it was shown in the recent 
very interesting paper
\cite{EG} by P. Etingof and V. Ginzburg. They have shown that for
the group $G$ of type $B_6$ and the multiplicity function $m$ equal 
to 0 on the long roots
and 1 on the short roots \cite{EG} $m$-harmonic polynomials
are linearly dependent over $S^G$. This raises the question of how 
exceptional this situation is.
Probably this never happens if the multiplicity $m$ is constant (in 
particular for all
Coxeter groups with only one conjugacy class of reflections).

Etingof and Ginzburg showed also that $Q_m$ is freely generated over $S^G$ by 
some homogeneous polynomials
which have the same degrees as $m$-harmonic polynomials. This means 
that one can use
our results to compute the Poincar\'e series for the algebra of quasiinvariants.
The fact that the
Poincar\'e polynomials $P(H_m,t)$ are palindromic corresponds to the 
Gorenstein property of $Q_m$
(see \cite{EG}).

It is interesting to compare the theory of $m$-harmonic polynomials 
with the classical case
$m=0$. Probably the main novelty is a more explicit role of the 
action of the group: as we have seen
each isotypical component of $H_m$ behaves differently when $m$ 
changes. In the classical case one
can write down an explicit formula for the Poincar\'e polynomial of $H_0$ without
any reference to representation theory (see Introduction) but for 
general $m$ this seems to be impossible.

We should mention that the problem of explicit description of 
$m$-harmonic polynomials is still
open. 
As we have shown, 
an essential step in its solution is  finding the 
solutions of the KZ equations \eqref{e-KZ0}
with values in the irreducible representations of $G$, 
but this is 
another interesting
open problem.

\section*{Acknowledgements}

We are grateful to Yu.~Berest, O. Chalykh, I. Cherednik, M. Shapiro and 
especially to M. Feigin
for useful discussions and
helpful comments. In particular M. Feigin attracted our attention to the papers
\cite{Kirillov},\cite{Solomon}
and I. Cherednik pointed out Opdam's paper \cite{Opdam2} to us.
We would like to thank also P. Etingof for comments
and information about his recent results with V. Ginzburg
\cite{EG}. About 
 Kerov's remarkable results
on the Plancherel measure of the symmetric groups $S_n$ we have 
learned from
G. Olshanski and A. Okounkov.

One of us (APV) is grateful to the Forschungsinstitut f\"ur Mathematik, ETH,
Zurich for the hospitality during April 2001 when this work has 
essentially been done.

This work has been partially supported by EPSRC (GR/M69548).

\end{document}